\documentstyle{article}
\begin{document}

\title{A cardinal invariant related to homogeneous families}
\author {Martin Goldstern \and Saharon Shelah}
\thanks{The first author was supported by an Erwin Schroedinger
  fellowship from the Austrian Science Foundation.}
\thanks{Publication E14.}

\date{Feb 23, 1997}

\long\def\hhh{\par\centerline{\vrule height 0.5pt width 0.5\hsize}\par}

\def\xx{{\bf ex}}
\def\cov{{\rm cov}}
\def\meager{{\cal M}}
\def\Y{{\bf Y}}
\def\B{{\cal B}}
\def\A{{\cal A}}
\def\F{{\cal F}}
\def\itm#1 {\item[{(#1)}]}
\def\dom{{\rm dom}}

\def\newx#1 {\newtheorem{#1}[thm]{#1}}
\newtheorem{thm}{thm}

\newx Lemma
\newx Fact
\newx Definition
\newx {Fact and Definition}
\newx {Definition and Fact}
\newx Remark
\newx Claim
\newx Conclusion

\maketitle

\begin{Definition}\label{def}
We call a family $\F \subseteq [\omega]^\omega $ independent if every 
boolean intersection 
$$ \bigcap_{A \in \F_0} A \cap \bigcap _{B \in \F_1} (\omega \setminus
B) \qquad\qquad
\mbox{ $\F_0 , \F_1 \subseteq \F$, disjoint, finite }
$$
is infinite. 

We call $\F$ homogeneous iff:  for every finite partial injective
functions $f: \omega \to \omega $ and $g: \F \to \F$  that satisfies 
$$ \forall  n\in \dom(f)\, \forall A\in \dom(g): \ 
  n \in A  \ \Leftrightarrow f(n) \in g(A)$$
there is a
permutation $\pi $ of $ \omega $ such that 
$\pi[A], \pi^{-1}[A]\in \F$ for all $A \in \F$, $\pi $ extends $f$ and
$\pi[A]=g(A)$ for all $A\in \dom(g)$. 

\end{Definition}

Explanation, motivation...

\begin{Definition}
Let $\xx= \min \{ |F|: $ $F$ is an independent family on $\omega$ that
cannot be extended to a homogeneous family.
\end{Definition}

It is easy to see that $\xx \ge \cov(\meager)$ (see below).  (This was
also pointed out by Brendle). 

We show here that in fact $\xx = \cov(\meager)$.  We will use the
following characterisation of $ \cov(\meager)$: 

\begin{Definition and Fact}
\begin{itemize} 
\item 
  If $A \subseteq \omega^\omega $, and $g\in \omega^ \omega $, we say
  that  $g$ 
  ``diagonalizes''
  $\F$ if $ \forall  f \in \F\, \exists^\infty n \ f(n)=g(n)$. 
\item By a theorem of Bartoszynski, $ \kappa < \cov(\meager)$ iff
  every family $\F \subseteq \omega ^ \omega $ of size $ \kappa $ can
  be diagonalized.  
\item If there is no $g$ which diagonalizes $\F$, $|\F|= \kappa $
 then we say that $\F$ ``witnesses'' $\cov(\meager) \le \kappa $.
\end{itemize}
\end{Definition and Fact}

\begin{Definition}
Define $ \kappa < \xx^*$ iff:  
 For any independent family $\A \subseteq [\omega]^\omega $ 
 of cardinality $\kappa$ there is a permutation $\pi:\omega \to \omega
 $ such that for all $A\in \A$ the set 
$$ \{ n \in A:  \pi(n)\in A \setminus \{n\} \}$$ is infinite. 
\end{Definition}

We will show
\begin{enumerate}
\itm 1 $\xx \le \xx^*$
\itm 2 $\xx^* \le \cov(\meager)$
\itm 3 $\cov(\meager)\le \xx$
\end{enumerate}

(1) and (3) are easy. 

\subsection*{Proof of (3)}

It is enough to show that: For any $\F$ of cardinality $<
\cov(\meager)$,  for any pair
$(f,g)$ as in definition \ref{def} there is an independent family
$\F' \supseteq \F$ and a permutation $\pi$ of $ \omega $ extending 
$f$ such that $\pi[A] = g(A)$.  

(We can then get a homogeneous family $\F* \supseteq \F$
 in $|\F|$ many steps, using a
bookkeping argument and taking care of a single demand $(f,g)$ in each
step.

Moreover, we may assume that $f$ and $g$ are in fact permutations (of
certain finite subsets of $ \omega $ and $\F$ respectively), and 
  $\F$ is sufficiently saturated:  For any finite disjoint sets $p, q
  \subseteq \omega $ there is $A \in \F$, $p \subseteq A \subseteq
  \omega \setminus q$.

Let $\hat \F = \F \setminus \dom(g)$, $\hat \omega = \omega \setminus
\dom(f)$.   

Note that the elements of $\dom(g)$ are generators of a free Boolean
algebra.  Extend $g$ to this algebra in a natural way ($g(A \cap B) =
g(A) \cap g(B)$, etc), and let 
$A_1, \ldots A_{n}$ be the atoms of this Boolean algebra.  
$g$ is a permutation of $\{A_1, \ldots A_{n}\}$.   We may assume that
$g$ acts onthe indices: $g(A_k) = A_{g(k)}$.

For each $k$ let $h_k : \omega \to A_k$ be an enumeration of $A_n
\setminus \dom(f)$.   

Now for any permutation  $c: \omega \to \omega $
define $\pi_c$ as follows: 
\begin{enumerate}
\item $\pi(i) = f(i)$ if $i \in \dom(f)$
\item $\pi(i) = h_{g(k)} ( \pi(h_k^{-1}(i)))$ if $i \in A_k\setminus
  \dom(f)$.  
\end{enumerate}

Note that $\pi_c$ is a permutation of $ \omega $, $\pi $ extends $f$,
and $\pi(A_k) = A_{g(k)} = g(A_k)$ for all $k\in \{1, \ldots, n\}$.  

Let $\F_c:= \{\pi_c^\ell[A]: A \in \F , \ell \in \omega \}$. 
  It remains to prove that we can choose $c$ such that $\F_c$ will be
  an independent family.    This follows easily from $|\F| <
  \cov(\meager)$.

\subsection*{Proof of (2)}



\begin{Definition}
Let $ X = \omega  \times \omega \times 2$.

Let $ Y = \{\rho_m: m\in \omega\}$ be the set of 
all partial finite functions from $X$ to $ \omega $. We say that $Y'
\subseteq Y $ is dense in $ Y $ iff 
$$ \forall \rho \in Y \,\, \exists \rho' \in  Y' : \rho \subseteq
\rho' $$

Let $ \kappa = \cov(\meager)$, witnessed by a family $(f_\alpha:
\alpha < \kappa)$.
\end{Definition}

 We will construct an independent family 
$(A_\alpha:  \alpha < \kappa)$ by induction, satisfying the following
condition: 
$$ (*)_{\beta}   \ \ \ 
\hbox{ If $B \in \B({\beta}  )$, then $\{\rho_n:n\in B\}$ is dense in
$Y$}$$
where $\B(\beta)$ is the collection of all finite 
boolean combinations of the
family $(A_\alpha: \alpha < \beta)$. 

Clearly $(*)_\beta$ is preserved in limit steps 
(i.e., $\forall \gamma < \beta \,(*)_\beta$ implies $(*)_\gamma$
whenever $\gamma$ is a limit ordinal). 

So let $\beta = \alpha +1$.   We have to define $A_\alpha$ such that
$(*)_\beta$ holds. 

To this end, consider the following forcing notion: 

\begin{Definition}
 $P_\alpha:= \{ w  \subseteq \omega $, $w $ finite, and: whenever
$m,n\in w$, $m < n$, then $\forall i\in \{0,1\}\, \exists k: 
\rho_n(m,k,i) = \eta_\alpha(m,k,i)$ $\}$. 
\end{Definition}

$P_\alpha$ is partially ordered by end extension.

\begin{Claim}
  For any $B \in \B(\alpha)$, any $\rho \in Y$ the sets
 $$D_{B,\rho} := \{ w \in  P_\alpha  : \exists n \in w \cap B, \rho
 \subseteq \rho_n \}$$
and 
 $$D'_{B,\rho} := \{ w\in P_\alpha : \exists n \in B \setminus w ,
         \rho \subseteq \rho_n \}$$
are dense in $P_\alpha $. 
\end{Claim}

Proof of the claim: Fix $B$, $\rho$, and let $u\in P_\alpha$. Wlog
$\min  B > \max u$. 
For each $m\in u$ and each $i\in \{0,1\}$ pick $k_{m,i} $ arbitrary
such that $( m, k_{m,i}, i) \notin \dom(\rho)$. 
Define  a finite partial function $\rho'$ on $X$ by requiring 
$\rho'(m,k_{m,i},i) =  \eta_\alpha (m,k,i)$ 
for all $m\in u$, $i\in \{0,1\}$. 

Now find $n \in B$ such that $\rho_n $ extends $\rho'$ and check that $u
\cup \{ n \} \in P_\alpha$. 

The density of $D'$ is similar (easier). 

\begin{Fact}  Assume that $G \subseteq G_\alpha $  is sufficiently generic,
i.e., $G$ meets all the dense sets mentioned above.  Then $G$ defines
a set $A_\alpha$ such that $(*)_{\alpha+1} $ holds. Moreover,
$$ (**)_\alpha:  \forall   i\in \{0,1\}\, \forall m<n: (m,n\in A)
\Rightarrow \exists k \rho_n(m,k,i)= f_\alpha(m,k,i)$$
\end{Fact}

\begin{Conclusion}
Let $\A = \{ A_\alpha: \alpha < \kappa \}$.  Then 
 there is no permutation such that $ \forall A \in \A $ the set $\{
n \in A: \pi(n) \in A \setminus \{n\}$ is infinite. 

Hence, $\A$ witnesses $ \xx^* \le \kappa $. 
\end{Conclusion}

Proof: 
Assume that $\pi$ is such a permutation.  Define $f:X \to \omega $ by 
$$\begin{array}{rcl}
f(m,k,0) &=& \rho_{\pi(m)}(m,k,0)\\
f(m,k,1) &=& \rho_{\pi^{-1}(m)}(m,k,1)\\
\end{array}
$$

We claim that now $f$ diagonalizes $\F$, which yields the desired
contradiction. Indeed, let $ \alpha <\kappa $. One of the two sets 
$$
\{ m \in A_\alpha : m < \pi(m) \in A_\alpha \} 
\ \ \ 
\{ n \in A_\alpha : n > \pi(n) \in A_\alpha \} 
$$
must be infinite. 

Case 1: $\{ m \in A_\alpha : m < \pi(m) \in A_\alpha \} $ is infinite.
For each $m$ in this set we can let $n:= \pi(m)$, (so 
$m<n$) and use $(**)_\alpha $ to show that 
$$  \exists k \ f(m,k,0) = \rho_n(m,k,0) = f_\alpha(m,k,0)$$

Case 2: $\{ n \in A_\alpha : n > \pi(n) \in A_\alpha \} $ is
infinite. Now for every $ n $ in this set we can let $m=\pi(n)$ (so
again $m<n$) and again use $(**)_\alpha $ to show that 
$$  \exists k \ f(m,k,1) = \rho_{\pi^{-1}(m)}(m,k,1) =
        \rho_n(m,k,1) =  f_\alpha(m,k,1)$$

Done.

\end{document}